\documentclass[a4paper,12pt]{amsart}
\usepackage{amssymb}
\usepackage{ifthen}

\nonstopmode \numberwithin{equation}{section}
\newtheorem{thm} {Theorem}
\newtheorem{cor} {Corollary}
\newtheorem{lem} {Lemma}
\newtheorem{prop} {Proposition}

\newtheorem{conj}{Conjecture}

\theoremstyle{definition}
\newtheorem{prob}{Problem}

\newenvironment{rem}{%
\bigskip
\noindent \textsl{{\sl Remark. }}}{\bigskip}
\newenvironment{rems}{%
\bigskip
\noindent \textsl{{\sl Remarks. }}}{\bigskip}

\newcounter {own}
\def\theown {\thesection       .\arabic{own}}

\newenvironment{pf}[1][]{%
 \vskip 3mm
 \noindent
 \ifthenelse{\equal{#1}{}}%
  {{\bf Proof. }}%
  {{\bf #1.} }%
 }%
{\qed\bigskip}

\newcounter{alphabet}
\newcounter{tmp}

\newcommand{\IR}{{\mathbb R}}
\newcommand{\ID}{{\mathbb D}}

\newcommand{\IC}{{\mathbb C}}

\def\be{\begin{equation}}
\def\ee{\end{equation}}

\newcommand{\bee}{\begin{enumerate}}
\newcommand{\eee}{\end{enumerate}}

\newcommand{\blem}{\begin{lem}}
\newcommand{\elem}{\end{lem}}
\newcommand{\bthm}{\begin{thm}}
\newcommand{\ethm}{\end{thm}}
\newcommand{\bcor}{\begin{cor}}
\newcommand{\ecor}{\end{cor}}
\newcommand{\beg}{\begin{examp}}
\newcommand{\eeg}{\end{examp}}
\newcommand{\begs}{\begin{examples}}
\newcommand{\eegs}{\end{examples}}
\newcommand{\bdefe}{\begin{defin}}
\newcommand{\edefe}{\end{defin}}
\newcommand{\bprob}{\begin{prob}}
\newcommand{\eprob}{\end{prob}}
\newcommand{\bei}{\begin{itemize}}
\newcommand{\eei}{\end{itemize}}

\newcommand{\bcon}{\begin{conj}}
\newcommand{\econ}{\end{conj}}
\newcommand{\bcons}{\begin{conjs}}
\newcommand{\econs}{\end{conjs}}
\newcommand{\bprop}{\begin{prop}}
\newcommand{\eprop}{\end{prop}}
\newcommand{\br}{\begin{rem}}
\newcommand{\er}{\end{rem}}
\newcommand{\brs}{\begin{rems}}
\newcommand{\ers}{\end{rems}}
\newcommand{\bo}{\begin{obser}}
\newcommand{\eo}{\end{obser}}
\newcommand{\bos}{\begin{obsers}}
\newcommand{\eos}{\end{obsers}}
\newcommand{\bpf}{\begin{pf}}
\newcommand{\epf}{\end{pf}}
\newcommand{\ba}{\begin{array}}
\newcommand{\ea}{\end{array}}
\newcommand{\beq}{\begin{eqnarray}}
\newcommand{\beqq}{\begin{eqnarray*}}
\newcommand{\eeq}{\end{eqnarray}}
\newcommand{\eeqq}{\end{eqnarray*}}

\newcommand{\ra}{\rightarrow}

\newcommand{\ds}{\displaystyle}

\begin{document}
\bibliographystyle{amsplain}
\title{On the Fekete-Szeg\"o problem for concave univalent functions}

\author{B. Bhowmik}
\address{B. Bhowmik, Department of Mathematics,
Indian Institute of Science, Bangalore-560012, India.}
\email{bappaditya.bhowmik@gmail.com}
\author{S. Ponnusamy}
\address{S. Ponnusamy, Department of Mathematics,
Indian Institute of Technology Madras, Chennai-600 036, India.}
\email{samy@iitm.ac.in}
\author{K.-J. Wirths}
\address{K.-J. Wirths, Institut f\"ur Analysis, TU Braunschweig,
38106 Braunschweig, Germany}
\email{kjwirths@tu-bs.de}

\date{ 
version: {\tt Nov 17, 2010}; File: $\tt FEK-Ar.tex$}
\begin{abstract}
We consider the Fekete-Szeg\"o problem with real parameter $\lambda$ for the class $Co(\alpha)$ of concave univalent functions.
\end{abstract}

\maketitle
\pagestyle{myheadings}
\markboth{B. Bhowmik, S. Ponnusamy and K.-J. Wirths
}{Fekete-Szeg$\ddot{\rm{o}}$ problem}

\section{Introduction}
Let ${\mathcal S}$ denote the class of all univalent (analytic) functions
\be\label{p9eq1}
f(z)=z+\sum_{n=2}^{\infty}a_n z^n
\ee
defined on the unit disk $\ID =\{z\in\IC:\,|z|<1\}$. Then the classical Fekete-Szeg\"{o} inequality, presented
by means of Loewner's method, for the coefficients of $f\in \mathcal{S}$ is  that
$$|a_3-\lambda a_2^2|\leq 1+2\exp (-2\lambda / (1-\lambda )) \quad \mbox{ for $\lambda\in [0,1)$}.
$$
As $\lambda \rightarrow 1-$, we have the elementary inequality $|a_3-a_2^2|\leq 1$. Moreover,
the coefficient functional
$$ \Lambda _\lambda (f)=a_3-\lambda a_2^2
$$
on the normalized analytic functions $f$ in the unit disk $\ID$ plays an important role in function
theory. For example, the quantity $a_3- a_2^2$
represents $S_f(0)/6$, where $S_f$ denotes the Schwarzian derivative
$(f''/f')'-(f''/f')^2/2$ of locally univalent functions $f$ in $\ID$.  In the literature, there exists a large
number of results about inequalities for $\Lambda _\lambda (f)$ corresponding to
various subclasses of $\mathcal{S}$.  The  problem  of maximizing the absolute value of the functional $\Lambda _\lambda (f)$
is called the Fekete-Szeg\"{o} problem.
In \cite{Koepf-87}, Koepf solved the Fekete-Szeg$\ddot{\rm{o}}$ problem for close-to-convex functions
and the largest real number $\lambda$ for which  $\Lambda _{\lambda} (f)$ is
maximized by the Koebe function $z/(1 - z)^2$ is  $\lambda= 1/3$, and later in
\cite{Koepf-87a} (see also \cite{London-93}),
this result was generalized for functions that are close-to-convex of  order $\beta$.
In  \cite{Pfluger-85}, Pfluger employed the variational method to give another treatment of the
Fekete-Szeg\"{o}  inequality which includes a description of the image domains under extremal functions.
Later, Pfluger \cite{Pfluger-86} used Jenkin's method to show that
$$|\Lambda _\lambda (f)| \leq 1+2|\exp (-2\lambda / (1-\lambda ))|, ~~ f\in {\mathcal S},
$$
holds for complex $\lambda$ such that ${\rm Re\,}(1/(1-\lambda))\geq 1$.
The inequality is sharp if and only if $\lambda$  is in a certain pear
shaped subregion of the disk given by
$$\lambda = 1-(u+itv)/(u^2 +v^2), ~~ −1\leq t\leq 1,
$$
where $u = 1- \log (\cos \phi)$ and $v = \tan \phi -\phi$ , $0< \phi< \pi/2$.

In this paper, we solve the Fekete-Szeg$\ddot{\rm{o}}$ problem
for functions in the class $Co(\alpha)$ of concave univalent functions, with real parameter $\lambda$.


\section{Preliminaries}

A function $f:\ID\to \IC$ is said to belong to the family $Co(\alpha)$  if
$f$ satisfies  the following conditions:
\bee
\item[(i)] $f$ is analytic in $\ID$ with the standard normalization $f(0)=f'(0)-1=0$.
In addition it satisfies $f(1)=\infty$.
\item[(ii)] $f$ maps $\ID$ conformally onto a set whose complement with respect to $\IC$ is
convex.
\item[(iii)]the opening angle of $f(\ID)$ at $\infty$ is less than or equal to
$\pi\alpha$, $\alpha\in (1,2]$.
\eee
This class has been extensively  studied in the recent years and for a detailed discussion
about concave functions, we refer to \cite{Avk-Wir-06, Avk-Wir-05,Pom-Cruz} and the references
therein.  We note that for  $f\in Co(\alpha)$, $\alpha\in (1,2]$, the closed set $\IC\backslash f(\ID)$
is convex and unbounded. Also, we observe that $Co(2)$ contains the classes $Co(\alpha)$,
$\alpha\in (1,2]$.


We recall the analytic characterization for functions in $Co(\alpha), \alpha \in (1,2]$: $f\in Co(\alpha)$ if and only if
${\rm\, Re}\, P_f(z)> 0$ in $\ID$,
where
$$P_f(z)=\frac{2}{\alpha-1}\,
 \left [\frac{(\alpha+1)}{2}\frac{1+z}{1-z}-1-z \frac{f''(z)}{f'(z)}\right ].
$$
In \cite{BPW-09}, we have used this characterization and proved the following theorem which
will be used to prove our result.

\bigskip
\noindent
{\bf Theorem A.} {\em
Let $\alpha\in (1, 2]$. A function $f\in Co(\alpha)$ if and only
if there exists a starlike function $\phi\in \mathcal{S}^{*}$ such that $f(z)=\Lambda _\phi (z)$,
where
$$\Lambda _\phi (z)
=\int_{0}^{z}\frac{1}{(1-t)^{\alpha+1}}\left(\frac{t}{\phi(t)}\right)^{(\alpha-1)/2}\, dt.
$$
}

We also recall a lemma due to Koepf \cite[Lemma 3]{Koepf-87}.

\medskip
\noindent
{\bf Lemma A.} {\em
Let $g(z)=z+b_2z+b_3z^2+\cdots \in \mathcal{S}^{*}$. Then $|b_3-\lambda b_2^2|\leq \max~\{1, |3-4\lambda|\}$
which is sharp for the Koebe function $k$ if $|\lambda-3/4|\geq 1/4$
and for $(k(z^2))^{1/2}=\frac{z}{1-z^2}$ if $|\lambda-\frac{3}{4}|\leq
1/4$.
}

Here $ \mathcal{S}^{*}$ denote the family of functions $g\in \mathcal{S} $ that map $\ID$ into domains that are
starlike with respect to the origin. Each  $ g\in\mathcal{S}^{*}$  is characterized by the condition
${\rm Re\,} (zg'(z)/g(z))>0$ in $\ID$. Ma and Minda \cite{Ma-Minda-92} presented the Fekete-Szeg\"o
problem for more general classes through
subordination, which includes the classes of starlike and convex functions, respectively.
In a recent paper, the authors in  \cite{CKS-07} obtained a new method of solving the
Fekete-Szeg\"o problem  for classes of close-to-convex functions defined in terms of subordination.

\section{Main Result and its Proof}
We recall from Theorem~A
 that $f\in Co(\alpha)$ if and only if there exists a function
$\phi (z)= z+ \sum_{n=2}^{\infty}\phi_n z^n\in \mathcal{S}^{*}$
such that
\be\label{p9eq2a}
f'(z)=\frac{1}{(1-z)^{\alpha+1}}\left(\frac{z}{\phi(z)}\right)^{\frac{\alpha-1}{2}},
\ee
where $f$ has the form given by (\ref{p9eq1}).
Comparing the coefficients of $z$ and $z^2$ on the both sides of the series expansion of (\ref{p9eq2a}), we obtain
that
\beqq
a_2&=& \frac{\alpha+1}{2}-\frac{\alpha-1}{4}\phi_2,\quad \mbox{and}\\
a_3 &=& \frac{(\alpha+1)(\alpha+2)}{6}-\frac{\alpha^2-1}{6}\phi_2-\frac{\alpha-1}{6}\phi_3+\frac{\alpha^2-1}{24}\phi_2^2,
\eeqq
respectively. A computation yields,
\beq\nonumber\label{p9eq5}
a_3-\lambda a_2^2 &= & \frac{(\alpha+1)^2}{4}\left ( \frac{2(\alpha+2)}{3(\alpha+1)}-\lambda \right )
+\frac{\alpha^2-1}{4}\left (\lambda-\frac{2}{3}\right )\phi_2\\
&& ~~~-\frac{\alpha-1}{6}\left [\phi_3-\left(\frac{2(\alpha+1)-3\lambda(\alpha-1)}{8}\right)\phi_2^2\right].
\eeq

\noindent
{\bf \underline{Case (1)}:} Let $\ds \lambda \in \left(-\infty, \frac{2(\alpha-3)}{3(\alpha-1)}\right]$.
We observe that the assumption on $\lambda$ is seen to be equivalent to
$$\frac{2(\alpha+1)-3\lambda(\alpha-1)}{8}\geq 1
$$
and the first term in the last expression is  nonnegative.
Hence, using Lemma~A for the last term in (\ref{p9eq5}), and noting that $|\phi_2|\leq 2$,
 we have from the equality (\ref{p9eq5}),
\beqq
|a_3-\lambda a_2^2|&\leq& \frac{(\alpha+1)^2}{4}\left ( \frac{2(\alpha+2)}{3(\alpha+1)}-\lambda \right )
+\frac{\alpha^2-1}{4}\left (\frac{2}{3}-\lambda \right ) \left|\phi_2\right|\\
&& ~~~ +\frac{\alpha-1}{6}\left| \phi_3-\left(\frac{2(\alpha+1)-3\lambda(\alpha-1)}{8}\right)\phi_2^2 \right|\\
&\leq & \frac{(\alpha+1)^2}{4}\left ( \frac{2(\alpha+2)}{3(\alpha+1)}-\lambda \right )
+\frac{\alpha^2-1}{2}\left (\frac{2}{3}-\lambda \right ) \\
&& ~~~ +\frac{\alpha-1}{6}\left(\frac{2(\alpha+1)-3\lambda(\alpha-1)}{2}-3\right).
\eeqq
Thus, simplifying the right hand expression gives
\be\label{neweq1}
|a_3-\lambda a_2^2| \leq \frac{2\alpha^2+1}{3}- \lambda \alpha^2, ~\mbox{ if }~ \lambda \in \left(-\infty, \frac{2(\alpha-3)}{3(\alpha-1)}\right].
\ee

\noindent
{\bf \underline{Case (2)}:} Let $\ds\lambda\geq \frac{2(\alpha+2)}{3(\alpha+1)}$
so that the first term in (\ref{p9eq5}) is  nonpositive.
The condition on $\lambda$ in particular gives $\lambda >2/3$ and therefore, our assumption on $\lambda$
implies that
$$ \frac{2(\alpha+1)-3\lambda(\alpha-1)}{8} < \frac{1}{2}.
$$
Again,  it follows from Lemma~A that
$$\left|\phi_3-\left(\frac{2(\alpha+1)-3\lambda(\alpha-1)}{8}\right)\phi_2^2\right|\leq 3-\frac{2(\alpha+1)-3\lambda(\alpha-1)}{2}.
$$
In view of these observations  and an use of the inequality $|\phi_2|\leq 2$, the equality (\ref{p9eq5}) gives
\beqq
|a_3-\lambda a_2^2| &\leq& -\frac{(\alpha+1)^2}{4}\left ( \frac{2(\alpha+2)}{3(\alpha+1)}-\lambda \right )
-\frac{\alpha^2-1}{2}\left (\frac{2}{3}-\lambda \right )\\
&& ~~~~ +\frac{\alpha-1}{6}\left(3-\frac{2(\alpha+1)-3\lambda(\alpha-1)}{2}\right).
\eeqq
Thus,  simplifying the right hand expression gives
\be\label{neweq2}
|a_3-\lambda a_2^2| \leq \lambda\alpha^2-\frac{2\alpha^2+1}{3}, ~\mbox{ if }~ \lambda\geq \frac{2(\alpha+2)}{3(\alpha+1)}.
\ee

The inequalities in both cases are sharp for the functions
$$f(z)= \frac{1}{2\alpha} \left [\left(\frac{1+z}{1-z}\right)^{\alpha}-1\right ].
$$

\noindent
{\bf \underline{Case (3)}:}
To get the complete solution of the Fekete-Szeg\"o problem,  we need to consider the case
\be\label{neweq3}
\lambda\in \left(\frac{2(\alpha -3)}{3(\alpha -1)},\,\frac{2(\alpha +2)}{3(\alpha +1)}\right).
\ee
Now, we deal with this case by using the formulas (\ref{p9eq2a}) and (\ref{p9eq5}) together with the representation formula
for $\phi\in \mathcal{S}^{*}$:
$$\frac{z\phi'(z)}{\phi(z)}=\frac{1+z\omega(z)}{1-z\omega(z)},
$$
where $\omega :\ID \ra \overline{\ID}$ is a function analytic in $\ID$ with the Taylor series
$$\omega (z)= \sum_{n=0}^{\infty}c_n z^n.
$$
Inserting the resulting formulas
$$\phi_2=2c_0~\mbox{ and }~\phi_3=c_1+3{c_0}^2
$$
into (\ref{p9eq5}) yields
\beqq
a_3-\lambda a_2^2&=& \frac{(\alpha+1)^2}{4}\left ( \frac{2(\alpha+2)}{3(\alpha+1)}-\lambda \right )
+\frac{\alpha^2-1}{2}\left (\lambda-\frac{2}{3}\right )c_0\\
&& ~~~ -\frac{\alpha-1}{6}\left[c_1+\frac{4-2\alpha+3\lambda(\alpha-1)}{2}{c_0}^2\right]\\
& =:& A + Bc_0+Cc_0^2+Dc_1,
\eeqq
where
$$\left\{
\ba{lll}
A &=&\ds \frac{(\alpha+1)(\alpha+2)}{6}-\lambda\frac{(\alpha+1)^2}{4},\\
B &=&\ds (\alpha^2-1)\left(\frac{\lambda}{2}-\frac{1}{3}\right),\\
C &=& \ds -\frac{(\alpha-1)\left(4-2\alpha+3\lambda(\alpha-1)\right)}{12},\\
D &=& \ds -\frac{\alpha-1}{6}.
\ea
\right.
$$
It is well-known that $|c_0|\leq 1$ and $|c_1|\leq 1-|c_0|^2$. Using this we obtain,
\beq\label{p9eq1a}
|a_3-\lambda a_2^2|&=&| A + Bc_0+Cc_0^2+Dc_1|\\\nonumber
&\leq&  |A + Bc_0+Cc_0^2|+ |D||c_1|\\\nonumber
&\leq&  |A + Bc_0+Cc_0^2|+ |D|(1-|c_0|^2).
\eeq

Let $c_0=re^{i\theta}$. First we search for the maximum of $|A + Bc_0+C{c_0}^2|$ where we fix $r$ and vary $\theta$.
To this end, we consider the expression

\vspace{8pt}
$|A + Bc_0+C{c_0}^2|^2$
\beqq
&=& |A+Bre^{i\theta}+Cr^2e^{2i\theta}|^2\\
&=&(A-Cr^2)^2+B^2r^2+(2ABr+2BCr^3)\cos\theta+4ACr^2(\cos\theta)^2\\
&=:& f(r,\theta).
\eeqq
Afterwards, we have to find the biggest value of the maximum function, if $r$ varies in the interval (0,1].

We need to deal with several subcases of (\ref{neweq3}).

\bigskip
\noindent
{\bf \underline{Case A}:} Let $\lambda\in \left(\frac{2(\alpha -3)}{3(\alpha -1)},\frac{2(\alpha -2)}{3(\alpha -1)} \right)$.
We observe that $C>0$, $B<0,$ and $A+Cr^2>0$ for $r\in [0,1]$. Hence the corresponding quadratic function
$$
h(x)=(A-Cr^2)^2+B^2r^2+2Br(A+Cr^2)x+4ACr^2x^2, ~x\in [-1,1],
$$
attains its maximum value for any $r\in (0,1]$ at $x=-1$. Therefore, our task is to find the maximum value of
$$
g(r)=A-Br+Cr^2+\frac{\alpha -1}{6}(1-r^2).
$$
The inequalities $g'(0)=-B$ and
$$
g'(1)=-B+2C-\frac{\alpha -1}{3}=\frac{\alpha -1}{6}(-6\lambda +4(\alpha -1))>0
$$
for $\lambda < \frac{2(\alpha -1)}{3\alpha}$ imply
$$
g(r)\leq g(1)=A-B+C=\frac{2\alpha^2 +1}{3}-\lambda \alpha^2.
$$

\bigskip
\noindent
{\bf \underline{Case B}:} If $\lambda =\frac{2(\alpha -2)}{3(\alpha -1)}$, then $C=0$ and $h$ is a linear function that has its maximum value at $x=-1$. The considerations of Case A apply and again we get the maximum value $g(1)$ as above.

\bigskip
\noindent
{\bf \underline{Case C}:} Let $\lambda \in \left(\frac{2(\alpha -2)}{3(\alpha -1)},\frac{2(\alpha -1)}{3\alpha} \right)$. Firstly, we prove that in this interval the quadratic function $h$ is monotonic decreasing for $x\in [-1,1]$. Since the function $h:\IR \ra \IR$ has its maximum at
$$
x(r)=\frac{-B(A+Cr^2)}{4ACr}=\frac{-B}{4}\left(\frac{1}{Cr}+\frac{r}{A}\right),
$$
it is sufficient to show that $x(r)$ is monotonic increasing and $x(1)<-1$.
The first assertion is trivial and the second one is equivalent to
$$j(\lambda)=\alpha^2(3\lambda-2)^2-4+3\lambda > 0
$$
for the parameters $\lambda$ in question. This inequality is easily proved.
Hence, we get the same upper bound as in Cases A and B.
In conclusion, Cases A, B and C give,
\be\label{neweq4}
|a_3-\lambda a_2^2| \leq \frac{2\alpha^2+1}{3}-\lambda\alpha^2, ~\mbox{ if }~
\lambda\in \left(\frac{2(\alpha -3)}{3(\alpha -1)}, \frac{2(\alpha-1)}{3\alpha} \right ).
\ee

\bigskip
\noindent
{\bf \underline{Case D}:} Let $\lambda \in \left[\frac{2(\alpha -1)}{3\alpha},\frac{2}{3} \right)$ and we may
factorize $j(\lambda)$ as
$$j(\lambda)=9\alpha^2(\lambda-\lambda _1)(\lambda-\lambda _2),
$$
where
\be\label{neweq8}
\lambda_1=\frac{4\alpha^2-1-\sqrt{8\alpha^2+1}}{6\alpha^2}~\mbox{ and }~\lambda_2=\frac{4\alpha^2-1+\sqrt{8\alpha^2+1}}{6\alpha^2}.
\ee
We observe that $\lambda_2>\lambda_1$.
For $\lambda \in  \left[\frac{2(\alpha -1)}{3\alpha},\lambda_1\right)$, the function $h$ has its maximum value at $x=-1$ and
the function $g$ has its maximum value at
$$
r_m=\frac{-B}{-2C+\frac{\alpha -1}{3}}\in (0,1].
$$
Hence, the maximum of the Fekete-Szeg\"o functional is
$$
g(r_m)=\frac{\alpha(10-9\lambda)-(3\lambda-2)}{9(2-\lambda )+3\alpha(3\lambda-2)}.
$$
 For $\lambda\in \left[\lambda_1,\frac{2}{3}\right)$, the number
$$
r_0=\frac{B}{2C\left(1+\sqrt{1-\frac{B^2}{4AC}}\right)}\in (0,1]
$$
is the unique solution of $x(r)=-1$ in the interval $(0,1]$. It is easily seen that $r_m<r_0$ for $\lambda < \frac{2}{3}$. Further,
$$
k(r)=\sqrt{h(x(r))}+\frac{\alpha -1}{6}(1-r^2)=(A-Cr^2)\sqrt{1-\frac{B^2}{4AC}}+\frac{\alpha -1}{6}(1-r^2)
$$
is monotonic decreasing for $r\geq r_0$. Hence, the maximum value of $|a_3-\lambda a_2^2|$ is $g(r_m)$ in this part of the interval in question, too. The extremal function maps $\ID$ onto a wedge shaped region with  an opening angle at infinity less than $\pi\alpha$ and one finite vertex as in Example 3.12 in \cite{BPW-08}.

\bigskip
\noindent
{\bf \underline{Case E}:} For $\lambda=\frac{2}{3}$, we have $B=0$ and $C=-\frac{\alpha -1}{6}$.
Hence the maximum is attained for $\cos\theta =0$ and any $r\in(0,1]$. In all these cases, we get
$$
|a_3-\lambda a_2^2|=\frac{\alpha}{3}
$$
as the sharp upper bound. The extremal functions map $\ID$ onto a region with an opening angle at infinity equal to $\pi\alpha$ and two finite vertices as in Example 3.13 in \cite{BPW-08}.

In conclusion, Cases D and E give,
\be\label{neweq5}
|a_3-\lambda a_2^2| \leq \frac{\alpha(10-9\lambda)-(3\lambda-2)}{9(2-\lambda )+3\alpha(3\lambda-2)}, ~\mbox{ if }~
\lambda\in \left(\frac{2(\alpha-1)}{3\alpha}, \frac{2}{3} \right ].
\ee

\bigskip
\noindent
{\bf \underline{Case F}:} Let $\lambda\in (\frac{2}{3},\lambda_2]$. Since $B>0$, the function
$x(r)$ is monotonic decreasing now. The number
$$
r_1=\frac{B}{-2C\left(1+\sqrt{1-\frac{B^2}{4AC}}\right)}\in (0,1]
$$
is the unique solution of the equation $x(r)=1$ lying in (0,1]. For $r< r_1$, we have $h(x)\leq h(1)$. We consider the function
$$
l(r)=A+Br+Cr^2+\frac{\alpha-1}{6}(1-r^2)
$$
and determine the maximum of this function to be attained at
$$
r_n=\frac{B}{-2C+\frac{\alpha-1}{3}}.
$$
It is easily proved that $r_n>r_1$. Since $k(r)$ is monotonic increasing, we get the maximum value of the Fekete-Szeg\"o functional in this case as
$$k(1)=(A-C)\sqrt{1-\frac{B^2}{4AC}}=\alpha(1-\lambda)\sqrt{\frac{12(1-\lambda)}{(4-3\lambda)^2-\alpha^2(3\lambda-2)^2}},
$$
which is attained for $c_0=e^{i\theta_0}$, where
$$
\cos\theta_0 = \frac{-B(A+C)}{4AC}.
$$
In this case, the extremal function $f$ is defined by the solution of the following complex differential equation
$$
f'(z)=\frac{(1-ze^{i\theta_0})^{\alpha-1}}{(1-z)^{\alpha+1}}.
$$
In conclusion, in this case, we have,
\be\label{neweq6}
|a_3-\lambda a_2^2| \leq\alpha(1-\lambda)\sqrt{\frac{12(1-\lambda)}{(4-3\lambda)^2-\alpha^2(3\lambda-2)^2}}, ~\mbox{ if }~
\lambda\in (2/3,\lambda_2].
\ee

\bigskip
\noindent
{\bf \underline{Case G}:} Let $\lambda \in \left(\lambda_2, \frac{2(\alpha +2)}{3(\alpha +1)}\right)$.
Since $x(1)<-1$ for these $\lambda$, the number
$$
r_2= \frac{B}{-2C\left(1-\sqrt{1-\frac{B^2}{4AC}}\right)}
$$
satisfies $x(r_2)=-1$ and $r_2\in (0,1)$.
For $r\leq r_2$, we can make similar considerations as in the preceding case, i.e. for $r\leq r_1$ the function
$l(r)$ takes the maximum value, and for $r\in (r_1,r_2]$, the function $k(r)$ plays this role.
For $r > r_2$, the point $x(r)$ does not lie  in the
interval $[-1,1]$. Hence, the maximum in question is attained for $x=-1$ or
$x=1$. We see that $A+C <0$ and $-A-Cr^2>0$ for the values of $\lambda$ that we are considering now,
the maximum of (\ref{p9eq1a}) is attained for $x=-1$, i.e. for
$c_0=-r$.
Hence, for $r \in (r_2,1]$ the maximum function is
$$
n(r)=-A+Br-Cr^2+\frac{\alpha -1}{6}(1-r^2).
$$
Since
$$
-C>\frac{\alpha -1}{6}~ \mbox{ and }~ B>0,
$$
we get $n(r)\leq n(1)$ in the interval in question and hence
$$|a_3-\lambda a_2^2|\leq n(1)=-A+B-C=\lambda \alpha^2-\frac{2\alpha^2 +1}{3}
$$
whenever $\lambda \in \left(\lambda_2, \frac{2(\alpha +2)}{3(\alpha +1)}\right)$.

Equations (\ref{neweq1}), (\ref{neweq2}),  (\ref{neweq4}),  (\ref{neweq5}), (\ref{neweq6}) and Case G give

\vspace{8pt}
\noindent
{\bf Theorem. }
{\it For $\alpha\in (1,2]$, let $f\in Co(\alpha)$ have the expansion $(\ref{p9eq1})$. Then, we have
$$\left|a_3-\lambda {a_2}^2\right|\leq
\left \{ \begin{array}{rl}
\ds \frac{2\alpha^2+1}{3}-\lambda\alpha^2 & \mbox{ for }\displaystyle \lambda \in
\left(-\infty,\frac{2(\alpha-1)}{3\alpha} \right ] \\
\ds \frac{\alpha(10-9\lambda)-(3\lambda-2)}{9(2-\lambda)+3\alpha(3\lambda-2)}
& \mbox{ for }\displaystyle \frac{2(\alpha-1)}{3\alpha}\leq \lambda \leq \frac{2}{3}\\
\ds \alpha(1-\lambda)\sqrt{\frac{12(1-\lambda)}{(4-3\lambda)^2-\alpha^2(3\lambda-2)^2}}
& \mbox{ for }\displaystyle \frac{2}{3}\leq \lambda \leq  \lambda_2\\
\ds \lambda\alpha^2-\frac{2\alpha^2+1}{3}
& \mbox{ for }\displaystyle \lambda\in \left [\lambda_2 ,\infty\right) ,
\end{array} \right .
$$
where $\lambda_2$ is given by {\rm (\ref{neweq8})}. To emphasize the fact that the bound is
a continuous function of $\lambda$ for any $\alpha$ we mention two different expressions for
the same bound for some values of $\lambda$. The inequalities are sharp.
}

\vspace{8pt}

The Fekete-Szeg\"o inequalities for functions in the class $Co(\alpha)$ for complex values
of $\lambda$ remain an open problem.

\vspace{8pt}
{\bf Acknowledgement:} The authors thank the referee for careful reading of the paper.

\end{document}